\documentclass{amsart}

\usepackage{amsmath,amssymb,amsfonts,amsthm,amsbsy,mathtools}
\usepackage[dvipsnames]{xcolor}
\usepackage{hyperref}
\usepackage{enumerate}
\usepackage[T1]{fontenc}
\usepackage{footnote,pifont,dsfont} 
\usepackage[mathscr]{eucal}
\usepackage{thmtools}
\usepackage{caption}
\usepackage{verbatim}
\usepackage{algorithm, algpseudocode}
\usepackage{listings}
\usepackage{fullpage}

\hypersetup{
    colorlinks=true,
    linkcolor=blue,
    citecolor=darkgray 
}

\newtheorem{theorem}{Theorem}

\theoremstyle{definition}
\newtheorem{defn}[theorem]{Definition}

\newtheorem{caution}[theorem]{CAUTION}

\theoremstyle{remark}
\newtheorem{rem}[theorem]{Remark}

\newtheorem{examplex}[theorem]{Example}
\newenvironment{example}
  {\pushQED{\qed}\examplex}
  {\popQED\endexamplex}

\definecolor{LightGrey}{RGB}{220, 220, 220}

\newcommand{\mC}{\ensuremath{\mathcal{C}}}

\newcommand{\mS}{\ensuremath{\mathcal{S}}}

\newcommand{\mV}{\ensuremath{\mathcal{V}}}

\newcommand{\C}{\mathbb{C}}
\newcommand{\Z}{\mathbb{Z}}
\newcommand{\R}{\mathbb{R}}
\newcommand{\N}{\mathbb{N}}
\newcommand{\Q}{\mathbb{Q}}

\newcommand{\bi}{\ensuremath{\mathbf{i}}}

\newcommand{\br}{\ensuremath{\mathbf{r}}}
\newcommand{\rr}{\ensuremath{\mathbf{r}}}

\newcommand{\vv}{\ensuremath{\mathbf{v}}}
\newcommand{\ww}{\ensuremath{\mathbf{w}}}

\newcommand{\bw}{\ensuremath{\mathbf{w}}}

\newcommand{\zz}{\ensuremath{\mathbf{z}}}
\newcommand{\one}{\ensuremath{\mathbf{1}}}
\newcommand{\zero}{\ensuremath{\mathbf{0}}}

\newcommand{\lgrad}{\ensuremath{{\nabla_{\log}}}}

\newcommand{\code}[1]{\texttt{\detokenize{#1}}}

\definecolor{codegreen}{rgb}{0,0.6,0}
\definecolor{codegray}{rgb}{0.5,0.5,0.5}
\definecolor{sageblue}{RGB}{94, 94, 255}
\definecolor{codepurple}{rgb}{0.58,0,0.82}
\definecolor{backcolour}{rgb}{0.95,0.95,0.92}

\lstdefinestyle{mystyle}{
    frame=tb,
    keywordstyle=\color{magenta},
    escapeinside={@}{@},
    backgroundcolor=\color{backcolour},   
    commentstyle=\color{codegreen},
    numberstyle=\tiny\color{codegray},
    stringstyle=\color{codegreen},
    basicstyle=\ttfamily\footnotesize,
    breakatwhitespace=false,         
    breaklines=true,                 
    captionpos=b,                    
    keepspaces=true,                 
    numbersep=5pt,                  
    showspaces=false,                
    showstringspaces=false,
    showtabs=false,                  
    tabsize=2,
    upquote=true,
}

\title[A SageMath Package for ACSV: Beyond Smoothness]{A SageMath Package for Analytic Combinatorics in Several Variables: Beyond the Smooth Case}


\author{Benjamin Hackl}
\address{Department of Mathematics and Scientific Computing, University of Graz, Austria}
\author{Andrew Luo}
\address{Cheriton School of Computer Science, University of Waterloo, Canada}
\author{Stephen Melczer}
\address{Department of Combinatorics and Optimization, University of Waterloo, Canada}
\author{Éric Schost}
\address{Cheriton School of Computer Science, University of Waterloo, Canada}

\usepackage[backend=bibtex]{biblatex}
\begin{document}

\begin{abstract}
The field of analytic combinatorics in several variables (ACSV) develops techniques to compute the asymptotic behaviour of multivariate sequences from analytic properties of their generating functions. When the generating function under consideration is rational, its set of singularities forms an algebraic variety -- called the \emph{singular variety} -- and asymptotic behaviour depends heavily on the geometry of the singular variety. By combining a recent algorithm for the Whitney stratification of algebraic varieties with methods from ACSV, we present the first software that rigorously computes asymptotics of sequences whose generating functions have non-smooth singular varieties (under other assumptions on local geometry). Our work is built on the existing \textsc{sage\_acsv} package for the SageMath computer algebra system, which previously gave asymptotics under a smoothness assumption. We also report on other improvements to the package, such as an efficient technique for determining higher order asymptotic expansions using Newton iteration, the ability to use more efficient backends for algebraic computations, and a method to compute so-called \emph{critical points} for \emph{any} multivariate rational function through Whitney stratification.
\end{abstract}

\maketitle

The field of \emph{analytic combinatorics in several variables (ACSV)}~\cite{Melczer2021,PemantleWilsonMelczer2024} develops techniques to study multivariate sequence using properies of their (multivariate) generating functions. In this document, our main goal is to take a rational function
\[ 
F(\zz) = \frac{G(\zz)}{H(\zz)} = \sum_{\bi\in\N^d}f_\bi \zz^\bi = \sum_{i_1,\dots,i_d\in\N^d}f_{i_1,\dots,i_d} z_1^{i_1}\cdots z_d^{i_d}
\]
defined by coprime polynomials $G,H\in\Z[\zz]$, and a \emph{direction vector} $\rr\in\N^d$ with non-zero coordinates, and determine asymptotics of the \emph{$\rr$-diagonal sequence} $\left(f_{n\br}\right)_{n\ge0}$.

\begin{caution}
The function $F$ is called \emph{combinatorial} if there exists a polynomial $P\in\Z[\zz]$ coprime with $H$ such that all but a finite number of the power series coefficients of $P(\zz)/H(\zz)$ are non-negative (for instance, this holds whenever all but a finite number of the $f_\bi$ are non-negative). For reasons described below, the main functionality of our package (automatic derivation of asymptotics) \textbf{assumes that $F$ is combinatorial}. It is an open problem whether combinatoriality is a decidable property, even in dimension one~\cite{OuaknineWorrell2014a}, so this is not verified by the package. All other conditions necessary to determine asymptotic behaviour are rigorously verified.
\end{caution}

\begin{rem}
Our results are more general than they may first appear. For instance, a diagonal of a $d$-variate algebraic function can be expressed as a modified diagonal of a $(d+1)$-variate rational function in an algorithmic way~\cite{GreenwoodMelczerRuzaWilson2022}, and the theory of ACSV implies that asymptotic behaviour of $f_\bi$ as $\bi\rightarrow\infty$ with $\bi/\|\bi\|\rightarrow\rr$ typically varies smoothly with~$\rr$, allowing for more general asymptotic results and even limit theorems~\cite{MelczerRuza2024} from the asymptotic behaviour of $\rr$-diagonals. Any multivariate generating function coming from an enumeration problem is necessarily combinatorial.
\end{rem}   

Unlike the univariate theory of analytic combinatorics~\cite{FlajoletSedgewick2009}, where many natural generating functions (including all rational, algebraic, and D-finite functions) admit a finite number of singularities, in dimension $d\geq2$ the \emph{singular variety} $\mV=\mV(H)$ defined by the vanishing of the denominator $H$ is a positive dimensional algebraic variety, and the geometry of $\mV$ plays a large role in the asymptotic analysis. The simplest case occurs when $H$ is square-free and $\mV$ is a smooth manifold, which happens precisely when $H$ and all of its partial derivatives don't simultaneously vanish. In 2023, the \code{sage_acsv} package~\cite{HacklLuoMelczerSeloverWong2023} was developed to rigorously determine leading asymptotics in this smooth setting. Here we report on new extensions that combine more advanced ACSV methods with a recent algorithm~\cite{HelmerNanda2022} for the Whitney stratification of algebraic varieties to give the first rigorous software computing asymptotics from multivariate generating functions with non-smooth singular varieties.

\begin{rem}
Although smoothness of $\mV(H)$ is a \emph{generic} property, meaning that it holds for all polynomials $H$ of a fixed degree except for those whose coefficients lie in a proper algebraic set, generating functions coming from combinatorial problems often exhibit non-generic behaviours (since having a combinatorial interpretation already makes them non-generic). We estimate that somewhere between a quarter and a half of generating functions described in the literature have non-smooth singular varieties.
\end{rem}

As before, the \code{sage_acsv} package can be installed in any recent SageMath installation
(preferably version \texttt{10.1} or later) by running the command 
\[ \code{sage -pip install sage-acsv} \]
from any Sage instance with access to the internet, or by downloading its source code from
\begin{center}
    \url{https://github.com/ACSVMath/sage_acsv}
\end{center}
Full documentation is available through Sage's built-in documentation system and at
\begin{center}
    \url{https://acsvmath.github.io/sage_acsv/}
\end{center}
This article discusses the functionality of version \code{v0.3.0}, released in April 2025.

\section{Computing Diagonal Asymptotics}

We begin with an example of a simple function whose singular variety is smooth.

\begin{example}
    The $(1,1)$-diagonal of the combinatorial rational function 
    \[F(x,y)=\frac{1}{1-x-y}=\sum_{a,b\ge0}\binom{a+b}{a}x^ay^b\] 
    forms the sequence $f_{n,n}=\binom{2n}{n}$. After installing our package, the code execution

\begin{lstlisting}[language=python]
sage: from sage_acsv import (get_expansion_terms,
....:     diagonal_asymptotics_combinatorial as diagonal)
sage: var('w x y z')
(w, x, y, z)
sage: diagonal(1/(1 - x - y))
1/sqrt(pi)*4^n*n^(-1/2) + O(4^n*n^(-3/2))
\end{lstlisting}
    verifies the required assumptions (other than combinatoriality) and proves that 
    \[\binom{2n}{n} = \frac{4^n}{\sqrt{\pi n}}\left(1+O\left(\frac{1}{n}\right)\right).\] 
    An optional argument to \code{diagonal_asymptotics_combinatorial} allows one to specify directions other than the default \emph{main diagonal} $\rr=\one$.
\begin{lstlisting}[language=python]
sage: diagonal(1/(1 - x - y), r=[2, 1])
0.866025403784439?/sqrt(pi)*(27/4)^n*n^(-1/2) + O((27/4)^n*n^(-3/2))
\end{lstlisting}
    The coefficient $\texttt{0.866025403784439?}$ is a printed numerical approximation of a constant that is computed exactly as an element of Sage's \code{Algebraic Field}. Asymptotics of the rational diagonal sequences we handle can be expressed as a sum of terms of the form $C\rho^n\pi^\alpha n^{\beta}$, where $C$ and $\rho$ are algebraic numbers and $\alpha,\beta\in\Q$. By default we return an element of \code{Asymptotic Ring}, but other options (such as a symbolic expression) are possible throught the optional \code{output_format} argument. New functionality added after~\cite{HacklLuoMelczerSeloverWong2023} also allows us to compute higher-order asymptotic expansions using the optional \code{expansion_precision} argument. Regardless of the chosen output format, or number of terms in the expansion, the \code{get_expansion_terms} function turns the expansion into a list of \code{Term} objects, from which the individual components can be more easily inspected.
\begin{lstlisting}[language=python]
sage: ex = diagonal(1/(1-x-y), r=[2, 1], expansion_precision=2); ex
0.866025403784439?/sqrt(pi)*(27/4)^n*n^(-1/2) - 0.08419691425682043?/sqrt(pi)*(27/4)^n*n^(-3/2) + O((27/4)^n*n^(-5/2))
sage: terms = get_expansion_terms(ex); terms
[Term(coefficient=0.866025403784439?, pi_factor=1/sqrt(pi), base=27/4, power=-1/2),
 Term(coefficient=-0.08419691425682043?, pi_factor=1/sqrt(pi), base=27/4, power=-3/2)]
sage: terms[0].coefficient.minpoly() # exact value of 0.86... is sqrt(3)/2
x^2 - 3/4
\end{lstlisting}
    Our higher-order asymptotic expansions use Newton iteration to implicitly compute high-order series expansions of analytic functions appearing in the analysis; here we can get 10 terms in the expansion in $\approx$1 second on a modern laptop.
\end{example}

We now provide some examples with non-smooth singular varieties.

\begin{example}
The coefficients of the generating function
\[ F(x,y) = \sum_{p,q\geq0}f_{p,q}x^py^q = \frac{1}{(1-x/3-2y/3)(1-2x/3-y/3)} \]
have an interpretation involving the number of winning choices in a single player game with a biased coin where the player must obtain $p$ heads and $q$ tails in $p+q$ flips~\cite[Example 12.25]{PemantleWilsonMelczer2024}. The derivation of $(r,s)$-diagonal asymptotics breaks down into different cases for 
\begin{itemize}
    \item $r/s < 1/2$, when the behaviour is dictated by a smooth point of the singular variety where the first denominator factor of $F$ vanishes but the second doesn't;
    \item $r/s > 2$, when the behaviour is dictated by a smooth point of the singular variety where the second denominator factor of $F$ vanishes but the first doesn't;
    \item $1/2 < r/s < 2$, when asymptotics are dictated by the \emph{non-smooth point} $(1,1)$ and we always have $f_{rn,sn}\sim 3$;
    \item and $r/s \in \{1/2,2\}$, when we have a \emph{non-generic direction}.
\end{itemize}
Our algorithm can automatically handle all cases except when $r/s \in \{1/2,2\}$ (see Definition~\ref{defn:contrib} below for more details on non-generic directions).

\begin{lstlisting}[language=python]
sage: F = 1/(1 - x/3 - 2*y/3)/(1 - 2*x/3 - y/3)
sage: diagonal(F, r=[1, 3])
6.531972647421808?/sqrt(pi)*(2048/2187)^n*n^(-1/2)
  + O((2048/2187)^n*n^(-3/2))
sage: diagonal(F, r=[4, 1])
3.952847075210475?/sqrt(pi)*(3125/3888)^n*n^(-1/2)
  + O((3125/3888)^n*n^(-3/2))
sage: diagonal(F, r=[1, 1])
3 + O(n^(-1))
\end{lstlisting}

Running \code{diagonal(F, r=[2,1])} throws the \code{ACSVException}: \code{Non-generic direction detected -} \\ \code{critical point [1, 1] is contained in 0-dimensional stratum.}
\end{example}

\begin{example}
\label{ex:latticepaths}
    The methods of ACSV have found great application in the enumeration of lattice path models restricted to convex cones~\cite[Chapters 4, 6, and 9]{Melczer2021}. For instance, in 2009 Bostan and Kauers~\cite{BostanKauers2009} guessed linear ODEs satisfied by the generating functions for 23 of the 79 non-isomorphic two-dimensional models restricted to a quadrant that take steps in $\{0,\pm1\}^2\setminus\{(0,0)\}$ and are not isomorphic to models restricted to halfspaces. From these linear ODEs, which were later rigorously derived~\cite{BostanChyzakHoeijKauersPech2017}, Bostan and Kauers combined the singularity analysis of \emph{D-finite functions} with heuristic guesses of numerically approximated constants to conjecture asymptotics of the 23 models. Due to the \emph{connection problem} for D-finite functions~\cite[Section 2.4.1]{Melczer2021}, several of these asymptotic conjectures stood until 2016, when they were resolved via ACSV~\cite{MelczerWilson2016}.

    As an explicit example, consider the sequence $s_n$ defined as the number of lattice paths that start at the origin, take $n$ steps in $\{(-1,-1),(1,-1),(0,1)\}$, and never leave the quadrant $\N^2$. The \emph{kernel method} for lattice path model enumeration~\cite[Chapter 4]{Melczer2021} implies that $s_n$ is the $(1,1,1)$-diagonal of the rational function
    \[ F(x,y,z) = \frac{(1 + x)(1-2zy^2(1+x^2))}{(1 - y)(1-z(x^2y^2 + y^2 + x))(1-zy^2(1+x^2))}, \]
    from which one can use the methods of \emph{creative telescoping} to compute~\cite[Model 15]{BostanChyzakHoeijKauersPech2017} that the generating function $S(t)$ of $s_n$ satisfies a linear ODE of order 5 with polynomial coefficients of maximum degree 16. Desingularizing this equation using the Sage \code{ore_algebra} package~\cite{KauersJaroschekJohansson2015,ChenKauersSinger2016} shows that $S(t)$ satisfies a linear ODE of order 43 whose leading coefficient is the polynomial
    \[ t^3 (t - 1/3)(t + 1)(t^2 - 1/8)(t^2 + 1/8), \]
    so the theory of D-finite function asymptotics~\cite[Section 2.4]{Melczer2021} then implies the existence of constants $\alpha,\beta,C,C_0,C_1,\dots$ such that
    \[ s_n = C \cdot 3^n n^\alpha \log^\beta(n) \sum_{k \geq 0} C_kn^{-k} + O((2\sqrt{2})^n). \]
    Using efficient algorithms~\cite{Mezzarobba2016} for the numeric analytic continuation of solutions of linear ODEs it is possible to rigorously approximate $C=0.000\dots$ to thousands of decimal places on a modern laptop, but it is unknown if it can be rigorously decided when such \emph{connection coefficients} are exactly 0 (which in this case would imply that $s_n$ has exponential growth at most $2\sqrt{2}$).

    Although we cannot proceed with univariate generating function methods\footnote{In this case one can use probabilistic arguments linking random walks to Brownian motion to prove that the exponential growth of $s_n$ is strictly less than 3, and then deduce it is $2\sqrt{2}$ via a generating function analysis, but this specialized argument for lattice paths cannot be automated for general sequences.}, with our package we can now rigorously compute asymptotics of $s_n$ directly from the multivariate function $F(x,y,z)$. Indeed, the code execution
\begin{lstlisting}[language=python]
sage: F = (1 + x)*(2*z*x^2*y^2 + 2*z*y^2 - 1)/((-1 + y)*(z*x^2*y^2 + z*y^2 + z*x - 1)*(z*x^2*y^2 + z*y^2 - 1))
sage: diagonal(F)
O(2.828427124746190?^n*n^(-2))
sage: diagonal(F, expansion_precision=2)
0.9705627484771406?/pi*(-2.828427124746190?)^n*n^(-2)
  + 32.97056274847714?/pi*2.828427124746190?^n*n^(-2)
  + O(2.828427124746190?^n*n^(-3))
\end{lstlisting}
    implies, after expressing the degree two algebraic numbers that appear in radicals, that
    \[ s_n = (12\sqrt{2} - 16)\frac{(-2 \sqrt{2})^n}{\pi n^2} 
            + (12\sqrt{2} + 16)\frac{(2 \sqrt{2})^n}{\pi n^2}. \]
    Note that we must set \code{expansion_precision=2} because the expected leading coefficient vanishes in this example, so by default we only conclude $s_n = O((2 \sqrt{2})^n/n^2)$.
\end{example}

\section{Theoretical Background and Other Functionality}

The methods of ACSV work by computing a (generically) finite collection of \emph{critical points} on the singular variety $\mV$, where local behaviour of $F$ \emph{could} impact asymptotics of the $\rr$-diagonal, and then refining them to points that actually do impact the asymptotics. The simplest situation is the square-free smooth case, when $H$ and its partial derivatives never simultaneously vanish, in which case the critical points in the direction $\rr$ are defined by the polynomial system 
\begin{equation}
    \begin{aligned}
        H(\bw) &= 0 \\
        r_kz_1H_{z_1}(\bw)-r_1z_kH_{z_k}(\bw) &= 0 \quad \text{for } 2 \le k \le d.
    \end{aligned}
\end{equation}
In general, the variety $\mV$ may not be smooth (due to self-intersections, cusp points, etc.) and must be partitioned into a finite collection of smooth sets.

\subsection*{Stratifications}
A \emph{Whitney stratification}~\cite{GoreskyMacPherson1988} of $\mV$ is a particular decomposition into smooth sets, called \emph{strata}, such that the local picture of $\mV$ near all points in the same stratum is consistent. Although originally described by Whitney in terms of limit behaviour of secant lines and tangent planes on $\mV$, an algebraic characterization of Whitney stratifications using \emph{conormal spaces} was given by L\^{e} and Teissier~\cite{TrangTeissier1988}. For our purposes, it is enough to note that there is a Whitney stratification defined by a collection of algebraic sets
\[ V(H) = X_d \supset X_{d-1} \supset \dots \supset X_0 \supset X_{-1} = \emptyset\] 
such that, for each $0\le k\le d$, the stratum $\mS_k = X_k\setminus X_{k-1}$ is either empty or a smooth $k$-dimensional manifold. We compute polynomial generators for the ideals $\mathcal{I}(X_k)$ using a recent algorithm of Helmer and Nanda~\cite{HelmerNanda2022}, in a function \code{whitney_stratification} that is used implicitly for ACSV purposes but may also be of independent interest.

\begin{example}
The following code computes a Whitney stratification of $\mV(y^2+z^3-x^2z^2)$.
\begin{lstlisting}[language=python]
sage: from sage_acsv import whitney_stratification
sage: R.<x, y, z> = PolynomialRing(QQ, 3)
sage: IX = Ideal(y^2 + z^3 - x^2*z^2)
sage: whitney_stratification(IX, R)
[Ideal (z, y, x) of Multivariate Polynomial Ring in x, y, z over Rational Field,
 Ideal (z, y) of Multivariate Polynomial Ring in x, y, z over Rational Field,
 Ideal (x^2*z^2 - z^3 - y^2) of Multivariate Polynomial Ring in x, y, z over Rational Field]
\end{lstlisting}
\end{example}

\begin{rem}
When $\mV$ is smooth then $\{\mV\}$ is its (trivial) Whitney stratification. 
\end{rem}

\begin{defn}
A collection of square-free polynomials $\{p_1,\dots,p_s\}$ is said to \emph{intersect transversely} if each of the varieties $\mV(p_k)$ is smooth (meaning each $p_k$ and its partial derivatives do not simultaneously vanish) and for each sub-collection $\{p_{\pi_1},\dots,p_{\pi_r}\}$ the gradients of the $p_{\pi_j}$ are linearly independent at any point of $\mV(p_{\pi_1})\cap\cdots\cap\mV(p_{\pi_r})$. If $H(\zz) = H_1(\zz)^{a_1}\cdots H_s(\zz)^{a_s}$ for positive integers $a_k$ and all collections of the irreducible polynomial factors $H_k$ intersect transversely then we say that $H$ has a \emph{transverse polynomial factorization}; if all $a_k=1$ then we also say that the factorization is \emph{square-free}.
\end{defn}

For ACSV purposes, when $H$ admits a transverse polynomial factorization then we can define \emph{flats} consisting of the simultaneous vanishing of all subsets of $\{H_1,\dots,H_s\}$ and make strata by taking each flat and removing all subflats. This process is more efficient than the general computations performed by \code{whitney_stratification}.

\begin{rem}
Our algorithm for diagonal asymptotics currently requires a transverse polynomial factorization, and can only compute higher-order asymptotic terms when $\mV$ is smooth or there is a square-free transverse factorization. Work is ongoing to compute higher-order terms in the general non-square-free case, and to reduce transverse factorization at the level of polynomials to transverse factorization in terms of local analytic functions near points determining asymptotics.
\end{rem}

\subsection*{Contributing Points}

The starting point of an asymptotic analysis is the Cauchy integral representation
\[ f_{n\rr} = \frac{1}{(2\pi i)^d} \int_{\mC} F(\zz) \frac{d\zz}{\zz^{n\rr+\one}}, \]
where $\mC = \{\zz\in\C^d : |z_1|=\cdots=|z_d|=\epsilon\}$ is a product of sufficiently small circles around the origin, and the function $\phi_\rr(\zz)=\zz^\rr$ captures the part of the Cauchy integrand that changes as $n\rightarrow\infty$. 

\begin{defn}
The \emph{critical points of $F$ on the stratum $\mS_k$} are defined to be the critical points of the map $\phi_\rr:\mS_k\rightarrow\C$ as a map of manifolds (i.e., the places where the differential of the restricted map $\phi_\rr{\big|}_{\mS_k}$ is rank deficient) and the set of \emph{critical points of $F$} is obtained by taking the union of the critical points on each stratum.
\end{defn}

If $(p_1,\dots,p_s)=\mathcal{I}(X_k)$ and $(q_1,\dots,q_r)=\mathcal{I}(X_{k+1})$ are the radical ideals defining the algebraic sets $X_k$ and $X_{k+1}$ then the critical points on the stratum $\mS_k$ are defined by the non-vanishing of $q_1,\dots,q_r$ together with the vanishing of $p_1,\dots,p_s$ and (as $\mS_k$ is a stratum of dimension $k$) the vanishing of the $(k+1)\times(k+1)$ minors of the matrix with rows $\nabla p_1, \dots, \nabla p_s, \nabla \phi$. Thus, critical points are always defined by a finite collection of polynomial equalities and inequalities. The function \code{critical_points} returns a list of the critical points of $F$ in the direction $\rr$ (with default $\rr=\one$) whenever this set is finite. In particular, this function \emph{does not assume that $F$ is combinatorial}. 

\begin{rem}
Because Sage's built-in functionality for ideal computations with Gröbner bases can be quite slow, we allow the user to change default settings using \code{ACSVSettings} to perform Gröbner computations using Macaulay2 (if installed on their system) and (often dramatically) speed up the computations. The path to Macaulay2 may also need to be specified, depending on the user's installation method and environment variables.
\end{rem}

\begin{example}
\label{ex:critpt}
The code execution
\begin{lstlisting}[language=python]
sage: from sage_acsv import ACSVSettings as AS, critical_points
sage: AS.set_default_groebner_backend(AS.Groebner.MACAULAY2)
sage: AS.set_macaulay2_path('/opt/homebrew/bin/M2')
sage: critical_points(1/(1 - (w + x + y + z) + 27*w*x*y*z))
[[1/3, 1/3, 1/3, 1/3],
 [-0.3333333333333334? + 0.4714045207910317?*I,
  -0.3333333333333334? + 0.4714045207910317?*I,
  -0.3333333333333334? + 0.4714045207910317?*I,
  -0.3333333333333334? + 0.4714045207910317?*I],
 [-0.3333333333333334? - 0.4714045207910317?*I,
  -0.3333333333333334? - 0.4714045207910317?*I,
  -0.3333333333333334? - 0.4714045207910317?*I,
  -0.3333333333333334? - 0.4714045207910317?*I]]
\end{lstlisting}
shows that $F(x,y,z,w) = 1/(1-(x+y+z+w)+27xyzw)$ admits 3 critical points in the direction $\one$: the point $(1/3,1/3,1/3,1/3)$, where $F$ has a so-called \emph{cone point singularity}, and a pair of points whose coordinates are complex conjugates, which are smooth points of the singular variety. Cone point singularities in even dimension at least 4 exhibit a \emph{lacuna phenomenon} that makes determining asymptotics particularly challenging~\cite{BaryshnikovMelczerPemantle2025}.
\end{example}

Although critical points are crucial to the analysis, not all critical points will affect asymptotic behaviour. We must therefore introduce some additional concepts.

\begin{defn}
    A \emph{minimal point} is an element $\ww\in\mV\cap\C_*^d$ with non-zero coordinates that is coordinate-wise minimal, meaning that there does not exist any $\vv \in \mV$ with $|v_j|<|w_j|$ for all $1 \leq j \leq d$.
\end{defn}

As in~\cite{HacklLuoMelczerSeloverWong2023}, we test for minimal critical points when $F$ is combinatorial using \emph{Kronecker representations} to reduce the necessary computations to bounding roots of univariate polynomials. We refer the interested reader to~\cite[Chapter 7]{Melczer2021} for theoretical details, or \cite[Section 2.2]{HacklLuoMelczerSeloverWong2023} for information on our implementation. Our test for minimality is what restricts us to studying combinatorial functions, as testing minimality for non-combinatorial functions is currently expensive enough to be non-feasible beyond very small cases. The function \code{minimal_critical_points_combinatorial} returns a list of the minimal critical points of $F$ in the direction $\rr$ (with default $\rr=\one$) whenever this set is finite.  

\begin{example}
The function $F_C(x,y,z,w) = 1/(1-(x+y+z+w)+Cxyzw)$ is known to be combinatorial whenever $C \leq 24$ (see~\cite{Yu2019}). When $C=24$ it can be shown using the \code{critical_points} function that $F$ admits 4 critical points in the direction $\one$, however the code execution
\begin{lstlisting}[language=python]
sage: from sage_acsv import minimal_critical_points_combinatorial
sage: minimal_critical_points_combinatorial(1/(1-(w+x+y+z)+24*w*x*y*z))
[[0.2961574755620697?, 0.2961574755620697?, 0.2961574755620697?, 0.2961574755620697?]]
\end{lstlisting}
shows that $F$ admits only one minimal critical point in this direction.
\end{example}

In the smooth case, if there are a finite number of minimal critical points then (under mild assumptions) they all have the same coordinate-wise modulus and determine dominant asymptotic behaviour. However, in the non-smooth case there are often different groups of minimal critical points with the same coordinate-wise moduli, with only the points in one group contributing asymptotically.  To determine the relevant points, we make use of an amazing fact about minimal points from the theory of hyperbolic polynomials. Given a polynomial $P$, let 
\[ (\lgrad P)(\zz) = \left(z_1P_{z_1}(\zz),\dots,z_dP_{z_d}(\zz)\right) \]
denote the \emph{logarithmic gradient of $P$}. If $\ww$ is a minimal point then $(\lgrad P)(\ww)$ is a scalar multiple of a real vector (see~\cite[Theorem 6.44]{PemantleWilsonMelczer2024}). 

\begin{defn}
\label{defn:contrib}
Suppose that $H$ admits the transverse polynomial factorization $H(\zz)=H_1(\zz)^{a_1}\cdots H_s(\zz)^{a_s}$ and let $\ww\in\mV\cap\C_*^d$ be a minimal critical point where (without loss of generality) precisely the factors $H_1,\dots,H_r$ vanish. For any $1 \leq k \leq d$ let $\vv_k \in \R$ be a real vector obtained by scaling $(\lgrad H_k)(\ww)$ by one of its non-zero entries. We say that $\ww \in \mV(H)$ is \emph{contributing} if there exist $\lambda_1,\dots,\lambda_r \geq 0$ such that
\begin{equation} \rr = \lambda_1 \vv_1 + \cdots + \lambda_r \vv_r. \label{eq:contrib_lambda} \end{equation}
The direction $\rr$ is called \emph{generic} if for any contributing point $\ww$ of $F$ the constants $\lambda_1,\dots,\lambda_r$ in~\eqref{eq:contrib_lambda} are all strictly positive.
\end{defn}

\begin{rem}
Our definition of contributing points does not depend on which non-zero entries of the $\lgrad H_k$ are used to define the $\vv_k$ in Definition~\ref{defn:contrib}. A critical point is a point where $\rr$ lies in the complex span of the $\lgrad H_k$, and a contributing point is a minimal critical point where $\rr$ lies in the non-negative real cone generated by the $\lgrad H_k$ after a suitable scaling. A non-generic direction $\rr$ is roughly one where a critical point of $\phi_\rr$ on a closure of a stratum lies in a stratum of lower dimension; as the name suggests, `most' directions are generic.
\end{rem}

The function \code{contributing_points_combinatorial} computes the contributing points in the direction $\rr$ (with default $\rr=\one$) when $H$ admits a transverse polynomial factorization; it fails with an \code{ACSVException} when a non-generic direction is detected.

\begin{example}
Returning to the lattice path model with steps $\{(-1,-1),(1,-1),(0,1)\}$ restricted to a quadrant, if $F$ is the multivariate function from Example~\ref{ex:latticepaths} then running
\begin{lstlisting}[language=python]
sage: minimal_critical_points_combinatorial(F)
[[1, 1, 1/3], [1, -0.7071067811865475?, 1/2], [1, 0.7071067811865475?, 1/2], [-1, 0.7071067811865475?*I, -1/2], [-1, -0.7071067811865475?*I, -1/2]]
sage: contributing_points_combinatorial(F)
[[1, 0.7071067811865475?, 1/2], [1, -0.7071067811865475?, 1/2], [-1, -0.7071067811865475?*I, -1/2], [-1, 0.7071067811865475?*I, -1/2]]
\end{lstlisting}
shows that $F$ has five minimal critical points -- the non-smooth point $(1,1,1/3)$ where two of the denominator factors vanish, and four smooth points -- with only the smooth points being contributing (smooth minimal critical points are always contributing points, when they exist). Note that only the two real points contribute to dominant asymptotics, but the other minimal critical points affect higher-order terms.
\end{example}

Once the contributing points are identified, \code{diagonal_asy} can compute the relevant asymptotic expansions under mild assumptions using~\cite[Theorems 9.1 and 9.2]{Melczer2021}. Given the combinatorial rational function $F(\zz)=G(\zz)/H(\zz)$ as input, the algorithm succeeds when: $H(\zero)\neq0$ so that $F(\zz)$ has a power series expansion at the origin; $H$ has a transverse polynomial factorization; $F$ admits a finite number of critical points, at least one of which is contributing; and the contributing points are \emph{non-degenerate}, which means that the Hessian matrices of certain implicitly defined analytic functions~\cite[Definition 9.10]{Melczer2021} are non-singular at the origin. All of these assumptions (except combinatoriality) are verified by the algorithm. If the numerator $H$ does not vanish for at least one contributing singularity then a sum of asymptotic terms will be returned. If $H$ vanishes at all such points then by default the algorithm will simply return a big-O bound; in this case, to compute dominant asymptotics a user should try using the \code{expansion_precision} argument to determine more terms in the asymptotic expansion.

\section{Acknowledgements}

AL partially supposed by an Ontario Graduate Scholarship, SM partially supported by NSERC Discovery Grant RGPIN-2021-02382, and ES partially supported by NSERC Discovery Grant RGPIN-2023-03463.

\printbibliography

\end{document}